\documentclass[journal]{IEEEtran}
\usepackage{amsmath,amssymb,amsfonts}
\usepackage{graphicx}
\usepackage{enumitem}
\usepackage{color}
\allowdisplaybreaks
\usepackage{titlesec}

\titlespacing\subsection{0pt}{8pt plus 2pt minus 2pt}{3pt plus 2pt minus 2pt}

\ifCLASSINFOpdf
 
\else
 
\fi

\begin{document}

\title{Uncertainty-aware Three-phase Optimal Power Flow based on Data-driven Convexification}

\author{Qifeng~Li,~\IEEEmembership{Senior Member,~IEEE}
\thanks{The author is with the Department
of Electrical and Computer Engineering, University of Central Florida, Orlando,
FL, 32816 USA.}
}


\maketitle
\begin{abstract}
This paper presents a novel optimization framework of formulating the three-phase optimal power flow that involves uncertainty. The proposed uncertainty-aware optimization (UaO) framework is: 1) a deterministic framework that is less complex than the existing optimization frameworks involving uncertainty, and 2) convex such that it admits polynomial-time algorithms and mature distributed optimization methods. To construct this UaO framework, a methodology of learning-aided uncertainty-aware modeling, with prediction errors of stochastic variables as the measurement of uncertainty, and a theory of data-driven convexification are proposed. Theoretically, the UaO framework is applicable for modeling general optimization problems under uncertainty.

\end{abstract}

\begin{IEEEkeywords}
Convex relaxation, data-driven, optimization under uncertainty, three-phase power flow.
\end{IEEEkeywords}
\IEEEpeerreviewmaketitle
\section{Introduction}
\IEEEPARstart{O}{PTIMIZATION} technologies have been widely used in many decision-making processes in the operation, control, and planning of power systems, such as optimal power flow (OPF). However, the increasing uncertainty introduced by distributed energy resources (DER) makes it extremely hard for operators to make accurate optimal decisions ahead of real time. There mainly exist three types of frameworks for modeling power system optimization problems that involve uncertainty: 1) stochastic framework, 2) robust framework, and 3) chance-constrained framework \cite{soroudi2013decision}. Unfortunately, these frameworks are rather computationally expensive for large-scale, highly-nonconvex problems. As a result, a large portion of existing works investigate proper assumptions to simplify these frameworks for power system applications \cite{aien2016comprehensive}. In contrast, based on regression analysis \cite{murphy2012machine}, this paper develops a novel uncertainty-aware optimization (UaO) framework using a new measurement of uncertainty that considers the prediction errors of stochastic variables (see Section III for more details).


Convex optimization \cite{boyd2004convex} has applications in a broad range of disciplines including power system engineering, mainly because: 1) many classes of convex optimization problems are computationally tractable as they admit polynomial-time algorithms; and, 2) it plays a fundamental role in the theories of both distributed optimization and bi-level optimization. The general idea is to relax the nonconvex objective functions and constraints into convex ones. However, a solution of the resulting convex problem may be infeasible to the original nonconvex problem due to the nature of relaxations, which now becomes one of the bottlenecks of this technology. In order to break these bottlenecks, this paper proposes a data-driven approach to construct convex relaxations with stronger tightness and lower complexity (see Subsection II-C for more details). The resulting convex relaxation is applied to convexify the developed UaO framework. The paper demonstrates the UoA framework on a three-phase optimal power flow (3$\phi$OPF) problem with uncertainty introduced by DERs and uncontrollable loads. The 3$\phi$OPF is balanced for transmission networks while unbalanced for distribution networks. It is worth noting that, theoretically, the proposed methods can be applied to general optimization problems under uncertainty.

\section{Theory of Data-Driven Convexification}

\subsection{Three-Phase Power Flow Equations}

In an OPF problem, the objective function is generally convex or linear. Thus, we focus on the main nonconvex constraints, i.e. the alternating current power flow (ACPF) equations, which are also considered as the mathematical model of power networks. Let $\mathcal{N}$, $\mathcal{E}$ and $\Phi$ denote the sets of buses, lines, and phases respectively. For each $i$, $j \in \mathcal{N}$, $ij \in \mathcal{E}$ and $\phi$, $\phi^\prime\in \Phi$, the formulation of three-phase ACPF (3$\phi$ACPF) equations \cite{hu2019ensemble} is given as
\begin{subequations} \label{PF}
\begin{align}
    &e_i^\phi\sum_{j\in \mathcal{N}}\sum_{\phi^\prime\in \Phi}(G_{ij}^{\phi \phi^\prime} e_j^{\phi^\prime}-B_{ij}^{\phi \phi^\prime} f_j^{\phi^\prime}) \nonumber \\
    &+f_i^\phi\sum_{j\in \mathcal{N}}\sum_{\phi^\prime\in \Phi}(B_{ij}^{\phi \phi^\prime} e_j^{\phi^\prime}
    +G_{ij}^{\phi \phi^\prime} f_j^{\phi^\prime})=p_{D,i}^\phi +p_{S,i}^\phi\\
     &f_i^\phi\sum_{j\in \mathcal{N}}\sum_{\phi^\prime\in \Phi}(G_{ij}^{\phi \phi^\prime} e_j^{\phi^\prime}-B_{ij}^{\phi \phi^\prime} f_j^{\phi^\prime}) \nonumber \\
    &-e_i^\phi\sum_{j\in \mathcal{N}}\sum_{\phi^\prime\in \Phi}(B_{ij}^{\phi \phi^\prime} e_j^{\phi^\prime}+G_{ij}^{\phi \phi^\prime} f_j^{\phi^\prime})=q_{D,i}^\phi 
    +q_{S,i}^\phi\\
    &e_i^\phi\sum_{\phi^\prime\in \Phi}[G_{ij}^{\phi \phi^\prime} (e_j^{\phi^\prime}-e_i^{\phi^\prime})-B_{ij}^{\phi \phi^\prime} (f_j^{\phi^\prime}-f_i^{\phi^\prime})]\nonumber \\
    &+f_i^\phi\sum_{\phi^\prime\in \Phi}[B_{ij}^{\phi \phi^\prime} (e_j^{\phi^\prime}-e_i^{\phi^\prime})
    +G_{ij}^{\phi \phi^\prime} (f_j^{\phi^\prime}-f_i^{\phi^\prime})]=p_{ij}^\phi \\
     &f_i^\phi\sum_{\phi^\prime\in \Phi}[G_{ij}^{\phi \phi^\prime} (e_j^{\phi^\prime}-e_i^{\phi^\prime})-B_{ij}^{\phi \phi^\prime} (f_j^{\phi^\prime}-f_i^{\phi^\prime})] \nonumber \\
    &-e_i^\phi\sum_{\phi^\prime\in \Phi}[B_{ij}^{\phi \phi^\prime} (e_j^{\phi^\prime}-e_i^{\phi^\prime})+G_{ij}^{\phi \phi^\prime} (f_j^{\phi^\prime}-f_i^{\phi^\prime})]=q_{ij}^\phi \\
    &(e_i^\phi)^2+(f_i^\phi)^2=v_i,
\end{align}
\end{subequations}
where $e$ and $f$ are the real and imaginary parts of voltage, while $p_{S,i}^\phi$ ($q_{S,i}^\phi$) and $p_{D,i}^\phi$ ($q_{D,i}^\phi$) denote the stochastic and deterministic components of active (reactive) power injections, respectively, at each bus. Let $|\mathcal{M}|=  (3\times|\mathcal{N}| +2\times|\mathcal{E}|)\times |\Phi|$, equations (\ref{PF}) can be compactly formulated as 
\begin{equation} \label{CompactPF}
     g_i(x)=x^{\rm{T}}A_i x=y_i=z_i+u_i, \;(i \in \mathcal{M})
\end{equation}
where the right-hand-side (RHS) vector $z=[\cdots,p_{D,i}^\phi,q_{D,i}^\phi,$ $\cdots,p_{ij}^\phi,q_{ij}^\phi,\cdots,v_i,\cdots]^ T$ denotes the RHS deterministic quantities while vector $u=[\cdots,p_{S,i}^\phi,q_{S,i}^\phi,\cdots,0,\cdots]^ T$ denotes the RHS stochastic quantities. The coefficient matrices $A_i$ in (1a)-(1d) are indefinite such that the quadratic functions in (\ref{CompactPF}) are nonconvex (except the ones corresponding to (1e)). Further define a set $\Omega = \{(x,\,y) |\,\underline{p}_i^\phi \le p_{D,i}^\phi + p_{S,i}^\phi \le \overline{p}_i^\phi,\,\underline{q}_i^\phi \le q_{D,i}^\phi +q_{S,i}^\phi \le \overline{q}_i^\phi, \, (p_{ij}^\phi)^2+(q_{ij}^\phi)^2 \le \overline{S}_{ij}, \, \text{and}\, \underline{v}_i \le v_i \le \overline{v}_i, \, \forall i \in \mathcal{N}\, \text{and}\, \phi \in \Phi \}$, then the feasible set of 3$\phi$ACPF (\ref{PF}) is $\Psi = \{(x,\,y)\in \Omega |\, (2),\,\forall i \in \mathcal{M}\}$. Note that $\Omega$ is convex while $\Psi$ is not. Moreover, the three-phase \textit{DistFlow} model \cite{baran1989optimal} of radial networks is also a nonconvex quadratic system that can be represented in the form of (\ref{CompactPF}). That means the proposed methods can be directly applied to \textit{DistFlow}-based 3$\phi$OPF.
\subsection{State-of-the-Art Convex Relaxations for AC Power Flow}
The ACPF can be formulated as quadratic models, like model (\ref{PF}), or trigonometric models. Most of the existing convex relaxations, e.g. the semidefinite (SD) \cite{bai2008semidefinite} and conic relaxations \cite{jabr2006radial}, are based on the quadratic ACPF models. Mature convex relaxations, such as the conic \cite{jabr2006radial} and convex hull \cite{li2017convex} relaxations, have been developed for the simplest scenario, i.e. the balanced ACPF of radial networks which can be formulated as the \textit{DistFlow} model \cite{baran1989optimal}. However, despite the satisfactory performance of these convex relaxations, radial (distribution) networks are barely balanced. For the rest scenarios, such as balanced ACPF in meshed networks and unbalanced ACPF (i.e. 3$\phi$ACPF) in radial networks, the most suitable one among the existing convex relaxations is the SD relaxation. However, the limitations and inexactness of the SD relaxation for ACPF have been widely observed \cite{lesieutre2011examining}. Moreover, the SD relaxation of large-scale OPF problems is computationally expensive. To mitigate the issues associated with the existing convex relaxations for ACPF, this paper explores to utilize the emerging machine learning technology in the process of constructing convex relaxations.
\subsection{Data-Driven Convex Relaxation}
In this subsection, a methodology of data-driven convex relaxation (DDCR) is established and applied to construct a tight convex quadratic relaxation of 3$\phi$ACPF model (\ref{CompactPF}) (i.e. (\ref{PF})). For simplicity, we start from a deterministic case, namely $u=0$. Let $\mathcal{D}$ denotes a training data set where the $k$th data point $D^{(k)}=( x^{(k)}, y^{(k)})$ denotes the real-time measurement of a historical operating point, we have $\mathcal{D} \subset \Psi$. The following regression algorithm is proposed to train  $\mathcal{D}$ to obtain a positive semi-definite (PSD) matrix $P_i$ and a complementary vector $B_i$ and scalar $c_i$ for each quadratic equations in (1a)-(1d):
\begin{subequations} \label{Regression}
\begin{align}
&\min_{P_i,B_i,c_i}\; \frac{1}{|\mathcal{D}|} \sum_kr^{(k)} \label{Reg1} \\
 \mathrm{s.t.}\; & (x^{(k)})^{\rm{T}}P_ix^{(k)}+B_i^{\rm{T}}x^{(k)}+c_i - y_i^{(k)} = m^{(k)} \le 0 \label{Reg2} \\
&\quad P_i \succeq 0 \label{Reg3} \\
&\quad \left[
\begin{array}{cc} 1     & m^{(k)} \\
 m^{(k)}  & r^{(k)} \\
\end{array}
\right]\succeq 0, \label{Reg4}
\end{align}
\end{subequations}
where $k=1,\,\ldots,\, |\mathcal{D}|$. Auxiliary variables $r^{(k)}$ are introduced for the purpose of formulating the regression model (\ref{Regression}) as a standard semidefinite programming (SDP) problem.
The dimensions of $P_i$, $B_i$, and $c_i$ are consistent with the dimensions of the corresponding quadratic equations in (1a)-(1d). Note that the $A$ matrix in (1e) is already PSD. Therefore, we don't need to train a $P$ for (1e). The optimization model (\ref{Regression}) is a standard SDP problem which can be effectively and globally solved by mature solvers like MOSEK, GUROBI, and CPLEX. 

Define a quadratic convex set:
\begin{equation} 
   \Theta =\{(x,\,y) \in \Omega\,|\, x^{\rm{T}}P_ix+B_i^{\rm{T}}x +c_i \le y_i, \forall i \in \mathcal{M} \}, \nonumber
\end{equation}
we have the following theorem.

\textbf{Theorem of Data-driven Convex Relaxation}. \textit{The set $\Theta$ is a convex relaxation of the feasible set $\Psi$ of the original three-phase AC power flow} (\ref{CompactPF}) \textit{if:}

\textit{a) the PSD matrices $P_i$, vectors $B_i$, and scalars $c_i$ ($i \in \mathcal{M}$) are obtained by training $\mathcal{D}$ using the regression algorithm (\ref{Regression}),}

\textit{b) $\mathcal{D}$ contains all extreme points\footnote{An extreme point of a convex set is a point in this set that does not line in any open line segment joining two points of this set \cite{bazaraa2013nonlinear}. We use this definition to define an extreme point of nonconvex sets.} of $\Psi$.}

\noindent
\textit{Proof}: Constraint (\ref{Reg2}) guarantees that each $D^{(k)} \in \mathcal{D}$ satisfies
\begin{equation}
    (x^{(k)})^{\rm{T}}P_ix^{(k)}+B_i^{\rm{T}}x^{(k)} +c_i \le y_i^{(k)}, \; i \in \mathcal{M} \nonumber
\end{equation}
which implies $\mathcal{D} \subset \Theta$. Therefore, $\Theta$ is a convex quadratic relaxation of $\mathcal{D}$ since $P_i$ ($i \in \mathcal{M}$) are PSD.

All extreme points of a feasible set are linearly independent according to the definition \cite{bazaraa2013nonlinear}. Suppose $\psi$ is an arbitrary point in $\Psi$, there must exist a vector of extreme points $X=[\theta_1,\theta_2,...,\theta_l]^{\rm{T}}$  of $\Psi$ and a vector of multipliers $\alpha=[\alpha_1, \alpha_2,...,\alpha_l]^{\rm{T}}$ that satisfy
\begin{equation}
    \psi=\alpha^{\rm{T}}X  , \nonumber
\end{equation}
where $0 \le \alpha_i \le 1$ ($i=1,\,\ldots,\, l$), $\sum_i^l \alpha_i =1$, and  $l$ equals to the dimension of the ($x, y$)-space. Since all $\theta_i \in \mathcal{D} \subset \Theta$ according to condition b), then $\psi \in \Theta$ due to the convexity of $\Theta$. Therefore, $\Psi \subset \Theta$ as $\psi$ is an arbitrary point in $\Psi$, which means $\Theta$ is convex relaxation of $\Psi$. \hfill$\square$

\noindent
\textbf{Remark 1}. Condition b) in the theorem of data-driven convex relaxation is not easy to strictly satisfy. However, under the concept of \textit{Big Data}, it is reasonable to assume that the data set $\mathcal{D}$ is big enough to represent the original feasible set $\Psi$, which implies $\Theta$ is highly close to a strictly convex relaxation of $\Psi$. Moreover, regression (\ref{Regression}) is a convex optimization problem that can be globally solved, which implies that $\Theta$ is the \textit{tightest} quadratic convex relaxation of $\Psi$.
\section{Learning-aided Uncertainty-aware Modeling}
This section introduces the idea and methodology of the proposed concept of learning-aided uncertainty-aware modeling (L-UaM) starting from the original 3$\phi$ACPF model (\ref{CompactPF}). For simplicity of the explanation, we ignore the controllable power injection $z$ which is deterministic and does not impact the analysis of uncertainty, namely $y=u$.
\subsection{Idea of L-UaM}
Suppose $\hat{u}$ is the real-time measurement of $u$ at time $t$, its system response is $\hat{x}$ as in equation (\ref{realtime}). Generally, the OPF is solved ahead of time $t$ (e.g. five minutes $\sim$ a day ahead) based on an ahead-of-real-time forecast $\tilde{u}$ as in (\ref{ahead}). One can consider that the uncertainty originates from the prediction error $\|\tilde{u}-\hat{u}\|$, since the resulting model output error $\|\tilde{x}-\hat{x}\|$ may lead to failures in OPF. The objective of this research is to construct an uncertainty-aware model $h$ that can approximate the actual state $\hat{x}$ at time $t$, which is a system response to real-time parameter $\hat{u}$, with only an ahead-of-real-time forecast $\tilde{u}$ as in (\ref{UaM}). For this purpose, a concept of L-UaM is proposed, of which the methodology is detailed in next subsection.
\begin{subequations}
\begin{align}
   & g_i(\hat{x})=\hat{u}_i \label{realtime} \\
   & g_i(\tilde{x})=\tilde{u}_i,\; (i \in \mathcal{M}) \label{ahead} \\
   & h_i(\hat{x}) \approx \tilde{u}_i \label{UaM}.
\end{align}
\end{subequations}

In the stochastic and chance-constrained frameworks, uncertainty is measured by probability distributions, while it is captured by a deterministic sets under the robust framework \cite{soroudi2013decision,aien2016comprehensive}. According to the above analysis, a ``good prediction" and a ``bad prediction" may have totally different impacts on the operation, control, and planning of power systems. However, the information of forecast is not considered in the measurement of uncertainty in these optimization frameworks, which increases the conservativeness. In this research, we use the prediction error $\|\tilde{u}-\hat{u}\|$ as the new measurement of uncertainty which will be incorporated into the following machine learning process. Researchers from the field of renewable generation and load forecast attempts to reduce the prediction errors $\|\tilde{u}-\hat{u}\|$, while this research focuses on reducing the model output errors $\|\tilde{x}-\hat{x}\|$ in the process of system modeling given the information of $\tilde{u}$.
\subsection{Methodology of L-UaM}
In order to learn an \textit{uncertainty-aware model} (UaM) of 3$\phi$ACPF, we first re-designed the historical training data set $\mathcal{D}$ to include the information of uncertainty. To avoid confusion, we denote the new data set as $\tilde{\mathcal{D}}$ of which the $k$th data point $\tilde{D}^{(k)}=(\hat{u}^{(k)},\hat{x}^{(k)},\tilde{u}^{(k)})$. Then, the following regression algorithm is designed to learn a hypothesis function $h$:
  \begin{equation} \label{URegression}
   \min_{h_i}\,\left \{ \frac{1}{|\tilde{\mathcal{D}}|} \sum_{\tilde{\mathcal{D}}}\left[h_i(\hat{x}^{(k)}) -\tilde{u}_i^{(k)}\right]^2 \right\}, \; (i \in \mathcal{M})
\end{equation} 
where the form of $h(\cdot)$ is pre-determined, and the cost function represents the mean of squared Euclidean distance between data points and the hypothesis function. The deterministic model $h(x)=u$ learned by regression (\ref{URegression}) is defined as a UaM of 3$\phi$ACPF (\ref{CompactPF}) (Note that the deterministic power component $z$ is ignored for simplicity in explanation. In the UaM training process of a actual power system, the term $z$ should be included). The UaM $h$ can capture the uncertainty of $u$ since it represents a mapping between $\tilde{u}$ and $\hat{x}$. Namely, it can produce a close approximation to the actual system state $\hat{x}$ at a future time $t$, with only an ahead-of-$t$ forecast values $\tilde{u}$ of the stochastic parameters $u$, as illustrated in equation (\ref{UaM}).

\noindent
\textbf{Remark 2}. Model $h$ resulting from (\ref{URegression}) is deterministic with ability of capturing the uncertainty of $u$ since it is obtained by training data set $\tilde{\mathcal{D}}$ that contains the historical forecast error $\|\tilde{u}-\hat{u}\|$ as the uncertainty measurement of $u$. Note that the forms of both cost function and the hypothesis function $h$ in regression (\ref{URegression}) are not unique, which can be customized as needed. For instance, to combine the L-UaM and DDCR processes, we will choose a convex quadratic hypothesis function as given in next section.

\section{Uncertainty-aware Convex Optimization}
\subsection{Uncertainty-aware Convex Modeling}
By incorporating the DDCR feature into the L-UaM process, this subsection aims at training a data-driven model of 3$\phi$ACPF that has two key features: convex and aware of uncertainty. For this purpose, like what we did in Section II, a hypothesis function $h$ of the following convex quadratic form is chosen:
\begin{equation} \label{UAM}
     h_i(x)=x^{\rm{T}}P_i x+B_i^{\rm{T}}x+c_i. \; (i \in \mathcal{M})
\end{equation}
Then, regression algorithms (\ref{Regression}) and (\ref{URegression}) are combined to produce the following regression algorithm
\begin{subequations} \label{U_regression}
\begin{align}
\min_{P,b,c}\; \text{(\ref{Reg1})} \quad  \mathrm{s.t.}\; &\text{(\ref{Reg3})},\;\text{(\ref{Reg4})} \; \text{and}  \\
& h_i(\hat{x}^{(k)}) - \hat{u}_i^{(k)} \le 0 \label{U_reg2} \\
& h_i(\hat{x}^{(k)}) - \tilde{u}_i^{(k)} = m^{(k)}, \label{U_reg3} 
\end{align}
\end{subequations}
to train the new data set $\tilde{\mathcal{D}}$. Constraint (\ref{U_reg3}) together with the objective function (\ref{Reg1}) implies that $h$ is an optimal mapping between real-time system response $\hat{x}$ and ahead-of-real-time forecast of stochastic parameter $\tilde{u}$. With the $P_i$, $B_i$ and $c_i$ inferred by regression (\ref{U_regression}), the epigraph (\ref{CUAM}) of $h$:
\begin{equation} \label{CUAM}
     \textbf{Epi}(h)=\{(x,y)\in \Omega\, |\,h_i(x) \le u_i, \forall i \in \mathcal{M} \}
\end{equation}
is a convex UaM of 3$\phi$ACPF (\ref{CompactPF}). Constraint (\ref{U_reg2}) indicates that convex set \textbf{Epi}($h$) encloses all data points in $\tilde{\mathcal{D}}$ which represent real-time operating points. It is worth pointing out that UaM $h$ uses the ahead-of-real-time forecast $\tilde{u}$ as input, which matches the fact that OPF is generally solved ahead-of-real-time.

\subsection{Uncertainty-aware 3$\phi$OPF}
With the typical objective that minimizes generation costs, the uncertainty-aware convex 3$\phi$OPF can be compactly formulated as
  \begin{equation} \label{OPF}
   \min_{z}\,\left \{ \sum_{i \in \mathcal{G}} \sum_{\phi \in \Phi}C(p_{D,i}^\phi):\,(x,z+\tilde{u}) \in \textbf{Epi}(h) \right\},
\end{equation}
where $\mathcal{G}$ denotes the generator set, and $P_i$, $B_i$, $c_i$ in \textbf{Epi}($h$) are inferred by training $\tilde{\mathcal{D}}$ using regression (\ref{U_regression}). A typical objective of OPF in distribution systems is to minimize the active power from transmission grids. It is worth noting that (\ref{OPF}) is a deterministic optimization problem which is much less-complex than the existing robust, stochastic, and chance-constrained optimization frameworks.
\section{Numerical Experiment}
\textit{1) Experiment design:} This section presents a numerical experiment on three scenarios of two balanced networks, i.e. a 5-bus and a 57-bus systems, and one unbalanced network, i.e. IEEE 34-bus distribution feeder. The three scenarios compared in this experiment are: 1) the original ACOPF with perfect predictions of stochastic power injections $u$; 2) the original ACOPF with inaccurate predictions $\tilde{u}$; and 3) the proposed convex UaO framework with the same inaccurate predictions $\tilde{u}$. Scenario 1 represents an ideal case which is considered as the reference. Scenario 2 simulates the current practice of solving OPF while scenario 3 is the proposed optimization framework for solving uncertainty-involved OPF. It is worth noting that this is a preliminary study based on simulated data since real-world data is currently not available. We will perform a more comprehensive study using real-world data in our future research.

\textit{2) Training of convex uncertainty-aware 3$\phi$ACPF models for the third scenario:} In the first step of generating the training data sets $\tilde{\mathcal{D}}$ for each test system, a set of power injection profiles $\hat{y}^{(k)}=\hat{z}^{(k)}+\hat{u}^{(k)}$ ($k=1,\cdots,50000$) are randomly produced. For each $\hat{u}^{(k)}$, there is a one-hour ahead forecast $\tilde{u}^{(k)}$ which is also randomly generated assuming a maximum forecast error of $\pm 30\%$, namely $|\tilde{u}^{(k)}-\hat{u}^{(k)}| \le 0.3|\hat{u}^{(k)}|$. Then, voltage profile $\hat{x}^{(k)}$ for each load profile of each system is obtained by solving PF (\ref{PF}). A convex uncertainty-aware 3$\phi$ACPF model is obtained for each of the three systems by training the corresponding data set $\tilde{\mathcal{D}}$ using regression (\ref{U_regression}).

\textit{3) Experiment results:} For each test system, the accuracy of scenario 3 is compared with scenario 2 considering scenario 1 as the reference on 50 load cases. Let $C_k^{(i)}$ denotes the optimal cost of the $k$th scenario in the $i$th load case, where $k=\{1,2,3\}$ and $i=1,\cdots,50$. The following average errors of objective values are calculated for scenarios 2 and 3 respectively
\begin{equation}
     E_k=\frac{1}{I}\sum_{i=1}^{I}\frac{|C_k^{(i)}-C_1^{(i)}|}{C_1^{(i)}}, \nonumber 
\end{equation}
where $k=\{2,3\}$ and $I=50$. The simulation results are tabulated in Table \ref{Results}.
\begin{table}[h]
\centering
\caption{Average Errors of Optimal Costs}
\label{Results}
\begin{tabular}{lccc}
\hline \hline
                       & \textbf{5-bus} & \textbf{34-bus} & \textbf{57-bus} \\ \hline
 $E_2$ & 20.21\% & 25.54\% & 21.43\% \\
$E_3$ & 3.72\% & 6.71\% & 1.66\% \\ \hline \hline
\end{tabular}
\end{table}

\textit{4) Analysis:} It can be observed from the above table that the proposed convex uncertainty-aware optimization framework (i.e. scenario 3) outperforms the original ACOPF (i.e. scenario 2) in terms of approximating the ideal case (i.e. scenario 1). The large errors of scenario 2 comes from the inaccurate forecast $\tilde{u}$ of the stochastic parameter $u$. Although the convex UaO framework also uses the inaccurate forecast $\tilde{u}$ as input, it can provide better solutions since the forecast errors are taken into account in the training process of the 3$\phi$ACPF constraints. The convex UaO framework is still not able to produce a strictly accurate solution. Nevertheless, it can be improved by training a larger, better data set according to its learning-based nature. Another important observation is that the proposed approaches are applicable to both balanced and unbalanced networks, and both radial and meshed networks.

Approximating a nonconvex optimization problem with its convex relaxations can significantly improve the computational efficiency while tightness is another important performance of a convex relaxation which is used to measure the accuracy of this approximation. The concept of convex UaO models is first proposed in this paper such that there not exist counterparts for comparison. Therefore, an independent numerical experiment in the deterministic environment is needed to evaluate the tightness of the proposed data-driven convex relaxation being compared with exiting deterministic convex relaxations. Due to page limit, we will seek to present such an evaluation in future publications based on real-world data.

\section{Conclusion and Future Work}
This paper develops a uncertainty-aware optimization (UaO) framework for modeling the three-phase optimal power flow (3$\phi$OPF) problem under uncertainty. The UaO is a convex, deterministic optimization framework that can capture uncertainty. A preliminary numerical experiment based on simulated data shows that the proposed framework has strong capability of mitigate the impacts of uncertainty on 3$\phi$OPF. In our future research, we will explore advanced machine learning technologies, such as ensemble learning \cite{murphy2012machine}, to improve the efficiency of UaO framework. Moreover, we will apply the UaO framework to model other power system optimzation problems other than OPF problems.





\ifCLASSOPTIONcaptionsoff
  \newpage
\fi



%
\bibliographystyle{unsrt}
\bibliography{main}

%








\end{document}